\documentclass[12pt,a4paper,twoside]{amsart}
\usepackage{graphics}
\usepackage{amsmath}
\usepackage{amssymb}

\bibstyle{plain}
\newtheorem{theorem}{Theorem}
\newtheorem{lemma}[theorem]{Lemma}
\theoremstyle{remark}
\newtheorem{remark}[theorem]{\bf Remark}

\newtheorem{example}{\bf Example}
\newcommand{\rng}{\{1,2,\ldots,\w,\mathfrak{c}\}}
\newcommand{\w}{\omega}
\newcommand{\co}{\mathfrak{c}}
\newcommand{\RR}{\mathcal{R}}
\newcommand{\sm}{\smallfrown}
\newcommand{\ra}{\vec{R}}
\newcommand{\Ra}{\vec{\RR}}
\begin {document}
\title{ Continuous functions taking every value a given number of times}
\title[Indicatrices of continuous functions]{ Continuous functions
 taking every value a given number of times}

\author{Aleksandra Kwiatkowska}
\address{Institute of Mathematics\\
 University of Wroclaw\\
Plac Grunwaldzki 2/4 \\50-384 Wroclaw\\Poland}
\email{ola\_171@wp.pl}

%\medskip
%\begin{center}
%{\bf Aleksandra Kwiatkowska}
%\end{center}

%%%%%%%%%%%%%%%%%%%%%%%%%%%%%%%%%%%%%%%%%%%%%%%%%%%%%%%%%%%%%%%%%%%%%%%%%%%%%%%%%%%%%%%%%%%%%

\begin{abstract}
We give necessary and sufficient conditions on 
a function $f\colon [0,1]\rightarrow \rng$ 
 under which  there exists a continuous function 
$F\colon [0,1]\rightarrow [0,1]$ such that for every $y\in[0,1]$,
 $|F^{-1}(y)|=f(y)$. 
\end{abstract}

\bigskip

\maketitle

\noindent {\bf 1. Introduction}

\vbadness=10000

  A function $f\colon [0,1]\rightarrow \rng$ will be called
a {\emph{ (Banach) indicatrix}} if there exists a continuous function 
$F\colon [0,1]\rightarrow [0,1]$ 
 such that $|F^{-1}(y)|=f(y)$ for each $y\in[0,1]$. 
We say that such a function $F$ is {\emph{described}} by $f$.
For each $y$, the set $F^y=\{x\colon F(x)=y\}$ will be called a \textit{section} of $F$.
In the paper we give a characterisation of indicatrices. 
Our result is a consequence   of the Darboux Property
and the
Sierpi\'nski-Mazurkiewicz Theorem, which in a special case  says that
for a continuous function $F\colon [0,1]\rightarrow [0,1]$, the set 
$F^{\co}=\{y\in[0,1]: |F^{-1}(y)|=\co\}$ is analytic.
We also give an explicit construction of a function described by any given indicatrix.
A general version of our result is given in  Theorem \ref{main2}.
In  Theorem \ref{main} we consider the case of  $F$ such that $F(0)=0$ and $F(1)=1$.

The term \textit{indicatrix} was introduced by Banach in \cite{Ba}.
He proved that the indicatrix of a continuous function is of Baire class 2 and
$V(F)=\int^1_0f(y)dy$, where $F$ is a coninuous function described by $f$
and $V(F)$ is the variation of $F$.

A construction of a continuous $F$
such that $F^{\co}$ is a given analytic set can be found in \cite{SM}. 
The case of $F^{\co}=[0,1]$ can be found in \cite{B} and \cite{G}.
Some partially results can also be found in \cite{C}.
A characterisation of indicatrices of Baire measurable and Lebesgue
measurable functions is in \cite {MN} and of Marczewski measurable ones is given in \cite{Ky}.
In \cite{KMM} there are characterized (under Analytic Determinancy) 
indicatrices of Borel measurable functions.

The paper is organized as follows. In Chapter 2 in Theorem \ref{main} 
we present a characterization of
indicatrices of continuous functions $F$ such that $F(0)=0$ and $F(1)=1$ and prove
that the conditions listed in this characterization are necessary. In  Chapter 3
we prove that they are also sufficient for $F$ with countable sections. In  Chapter 4
we deal with the general case and prove that the conditions listed in Theorem \ref{main}
are sufficient for $F$ which may have uncountable sections. In  Chapter 5 we show
that the assumption $F(0)=0$ and $F(1)=1$ may be dropped.

\bigskip\bigskip 
%%%%%%%%%%%%%%%%%%%%%%%%%%%%%%%%%%%%%%%%%%%%%%%%%%%%%%%%%%%%%%%%%%%%%%%%%%%
\noindent {\bf 2. The characterisation of indicatrices}

Define  for $\kappa\in\rng$
$f^{\kappa}=\{y\in[0,1]\colon f(y)=\kappa\}$, 
$f^{\geq\kappa}=\{y\in[0,1]\colon f(y)\geq\kappa\}$, 
and likewise define $f^{<\kappa}$, $f^{>\kappa}$, etc.

\begin{theorem}\label{main}
A function  $f\colon [0,1]\rightarrow \rng$  is an indicatrix of a continuous function 
$F\colon[0,1]\rightarrow[0,1]$
 such that $F(0)=0$ and $F(1)=1$
if and only if 

\noindent 

\emph{(1)} For all $y\in (0,1)$ such that $f(y)\in\w$,

 \ there exists $\varepsilon>0$ such that $y-\varepsilon<y_1<y<y_2<y+\varepsilon$
implies
$$ (*)  \qquad f(y_1)+f(y_2)\geq 2f(y),$$
and the equality can hold only if at least one (equivalently both)
of $f(y_1)$ and $f(y_2)$ is odd.

\indent Moreover, there exists $\varepsilon>0$ such that $0<y_2<\varepsilon$
implies $1+f(y_2)\geq 2f(y)$ and  $1-\varepsilon<y_1<1$
implies $f(y_1)+1\geq 2f(y)$.

\noindent \emph{(2)} 
For all $y\in(0,1)$ such that $f(y)=\w$ and for all but countably many $y$  such that 
$f(y)=\co$,

\ for every  $n\in\w$ there exists $\varepsilon>0$ such that 
$y-\varepsilon<y_1<y<y_2<y+\varepsilon$ implies
$$(**) \qquad  f(y_1)+f(y_2)\geq n.$$ 

\indent Moreover, for evry $n$ there exists $\varepsilon>0$ such that $0<y_2<\varepsilon$
implies $1+f(y_2)\geq 2f(y)$ and  $1-\varepsilon<y_1<1$
implies $f(y_1)+1\geq 2f(y)$.

\noindent \emph{(3)} The set $f^{\co}$ is analytic.

\end{theorem}

\begin{proof}[Proof of necessity]

 (1) Assume that $0<y<1$ and $f(y)=n$. Let 
$a_1<a_2<\ldots <a_n$ be all points which $F$ sends to $y$.
Consider intervals $[a_k, a_{k+1}], k=1,2,\ldots,n-1$ and choose
$b_k\in(a_k,a_{k+1})$. 
Let $\varepsilon=\min\{|F(b_k)-y|\colon k=1,2,\ldots,n-1\}$ and fix $y_1\in(y-\varepsilon,y)$
and $y_2\in(y,y+\varepsilon)$. It follows from the Darboux
Property that in every interval $(0,a_1), (a_k, b_k), (b_k, a_{k+1}), (a_n,1)$ either
 $y_1$ or $y_2$ is assumed at least once. There are $2n$ intervals, therefore we get
the  inequality in (*).

Now suppose that we have the equality in (*). Then in the intervals $(a_k, a_{k+1}),
 k=1,2,\ldots,n$,
exactly one of the values $y_1$ and $y_2$ must be taken exactly twice, and in the intervals 
$(0,a_1), (a_n,1)$ exactly once. If $f(y_1)$ and $f(y_2)$ were both even, then either $y_1$ 
or $y_2$
would be assumed in both of the intervals $(0,a_1), (a_n,1)$.
Then $F(0)\neq 0$ or $F(1)\neq 1$.
 We get a contradiction.

We similarly deal with the remaining part of (1).

(2) Assume  that $0<y<1$, $f(y)\in\{\w,\co\}$ and $F$ is not constant and equal $y$ in any
interval. Fix $n\in\w$. Choose $n$ points 
$a_1<a_2<\ldots <a_n$ which $F$ sends to $y$.
 $F$ is not constant on any interval
$[a_i,a_{i+1}]$. Similarly as before we get
$\varepsilon>0$ such that $y-\varepsilon<y_1<y<y_2<y+\varepsilon$
implies $f(y_1)+f(y_2)\geq 2n$, which gives (**). 

We similarly deal with cases $y=0$ and $y=1$. 

(3) This is a special case of the Sierpi\'nski-Mazurkiewicz Theorem (a proof can be found in
\cite{S} in subchapter 4.3, in \cite{K} or in \cite{SM}).
\end{proof}

Note that in the proof of (1) and (2) it suffices that $F$ is a Darboux fuction.
Therefore (1) and (2) characterise indicatices of Darboux functions $F$ satisfying 
$F^{\co}=\emptyset$, $F(0)=0$ and $F(1)=1$.

\medskip

Now we focus on the proof of sufficiency.

A function $f$ which satisfies conditions (1),(3) and satisfies (**) of (2) of
Theorem \ref{main} for all $y\in f^{\geq \w}$ will be called 
a {\emph{ preindicatrix}}. 

%%%%REMARKS%%%%%%%%%%%%%%%%%%%%%%%%%%%%%%%%%%%%%%%%%%%%%%%%%%%%%%%%%%%%%%%%

\begin{remark}
Without loss of generality
we may drop the restriction "but countably many" from condition (2) of Theorem \ref{main}
and consider only preindicatrices.

Indeed, let $S$ be the countable set of $y\in[0,1]$
such that $f(y)=\co$ and $f$ is bounded in some surrounding of $y$.   

Put for $y\in(0,1]$
$$f(y^-)=\liminf_{z\rightarrow y^-}f(z)$$ 
and for $y\in[0,1)$
$$f(y^+)=\liminf_{z\rightarrow y^+}f(z),$$ 
Let also $f(0^-)=1$ and  $f(1^+)=1$ 

Note that $f(y^+),f(y^-)\leq\w$.

For $y\in S$ put $f'(y)=\min\{f(y^-),f(y^+)\}-1$ if $f(y^+)$ and $f(y^-)$ are 
both even and equal and put
$f'(y)=\min\{f(y^-),f(y^+)\}$ otherwise. Let $f'(y)=f(y)$ outside $S$. Then
$f'$ is a preindicatrix.

Suppose we can find 
$F'$  described by $f'$ such that $F'(0)=0$ and $F'(1)=1$. 
At each level $y\in S$ inject into the graph of $F'$ a small horizontal interval
and shrink the resulting graph so that it has domain equal to [0,1]. In this way
we get a function $F$ such that
$\left|F^{-1}(y)\right|=\co$
for $y\in S$ and $\left|F^{-1}(y)\right|=\left|F'^{-1}(y)\right|$ outside $S$.
Therefore $F$ is described by $f$.
\end{remark}

A preindicatrix whose range is contained in $\{1,2,3\}$
will be called a  {\emph{ simple preindicatrix} }.
For a simple preindicatrix $p$ any maximal interval of $p^3$ will be called
\textit{interval of type three}.

\begin {remark}\label{three}
A function $p$ is a simple preindicatrix iff $p\colon[0,1]\rightarrow \{1,2,3\}$ and

a) $p^3$ is open; 

b) if $p(y)=2$ then $y$ is an endpoint of some interval of type three.

\end{remark}

\begin{lemma}\label{series}
Assume that $f$ is a preindicatrix. Then there exist a preindicatrix $f'$
and a simple preindicatrix $p$ such that for all $y\in [0,1]$,
$$f(y)=(p(y)-1)+f'(y).$$
(in the above equality if $f(y)\geq\w$ put $f(y)=\infty$).

\end{lemma}

\begin{proof}
Set $p\colon[0,1]\rightarrow\{1,2,3\}$ as follows
\begin{displaymath}
p(y)=\left\{
\begin{array}{ll}
1 & \textrm{if}\quad f(y)=1\\
3 & \textrm{if} \quad\exists_{\varepsilon>0}\forall_{z\in(y-\varepsilon, y+\varepsilon)}f(z)>2\\
2 & \textrm{otherwise},
\end{array}\right. 
\end{displaymath}

 and let
$$ f'(y)=\left\{
\begin{array}{ll}
f(y)-p(y)+1 &  \textrm{if}\quad f(y)<\w\\
\infty & \textrm{if}\quad f(x)\geq\w .
\end{array}\right. $$

\medskip
We show that the functions $p$ and $f'$ are preindicatrices.

First we check a) and b) of Remark \ref{three} for $p$.

Fix $y\in[0,1]$. Assume first that $p(y)=3$. Then for some neighbourhood 
$(y-\varepsilon, y+\varepsilon)$ of $y$ we have 
 $f\geq 3$. In the same interval $p\equiv3$. This implies that the set 
$\{y\in[0,1]\colon p(y)=3\}$ is open.

If $p(y)=2$ then $f(y)\geq 2$. Since $f$ is a preindicatrix,
for some $\varepsilon>0$ $f\geq 3$ in $(y-\varepsilon, y)$ or in $(y, y+\varepsilon)$. Then $p\equiv3$
in this neighbourhood, as desired.

\medskip
\noindent Now we prove that $f'$ is a preindicatrix. Fix $y\in(0,1)$.

\noindent Let first $p(y)=2$ and $f(y)<\w$. 
We apply (1) of the already proved part of 
Theorem \ref{main} to $f$ and fixed $y$ and take the resulting $\varepsilon>0$. 
If $f(y)=2$, then $f'(y)=1$ and the required inequalities holds. 
Otherwise, without loss of generality there is $y_0\in(y-\varepsilon, y)$ such  that 
$f(y_0)\in\{1,2\}$.  
Then if $y_2\in(y, y+\varepsilon)$ we obtain $$2f(y)\leq f(y_2)+f(y_0)$$
and equality can hold only if $f(y_2)$ and $f(y_0)$ are both odd, in particular when 
$f(y_0)=1$. Hence we get $$2f(y)\leq 1+f(y_2).$$
Therefore for $y_1\in(y-\varepsilon, y)$ and $y_2\in(y, y+\varepsilon)$ we get
$$\begin{array}{l}
f'(y_1)+f'(y_2)\geq 1+\left((1+f(y_2))-p(y_2)\right)\geq  1+2f(y)-p(y_2)\geq \\
1+2f(y)-3\geq 2(f(y)-p(y)+1)\geq 2f'(y).\end{array}$$

If $f'(y_1)$ and $f'(y_2)$ are both even, then the first of above inequalities is
 strict,
then $$f'(y_1)+f'(y_2)>2f'(y).$$

If $p(y)=1$ then $f(y)=1$.  Thus also
$f'(y)=1$. Hence for this $y$ the demanded inequalities hold.

If $p(y)=3$ and $f(y)\in\w$ then there exists a neighbourhood $(y-\varepsilon, y+\varepsilon)$
of $y$
in which $p$ receives value three. Therefore for $x\in(y-\varepsilon, y+\varepsilon)$ we have
$f'(x)+2=f(x)$. Hence the demanded inequalities follow from the similar ones for $f$.

The remaining case, i.e. when $y\in\{0,1\}$ is much easier to check (we omit it).

Finally, if $f(y)=\infty$, then also $f'(y)=\infty$. We apply (2) of the already proved part of 
Theorem \ref{main} to $n+2$ 
and $f$ and get the required inequality for $f'$.
This completes the proof of the Lemma.
\end{proof}   

\medskip

The following theorem will allow us to consider in the construction a sequence of
simple preindicatrices instead of a preindicatrix.

\begin{theorem}\label{series2}
Assume that $f$ is a preindicatrix.
Then there exist simple preindicatrices $p_0\geq p_1\geq{}
p_2\geq\ldots$ such that   for all $y\in[0,1]$ 
\begin{equation}\label{eq:iser}
f(y)=1+\sum_{i=0}^{\infty}(p_i(y)-1) ,  \qquad \textrm{ if } f(y)\in\w;
\end{equation}
\begin{equation}\label{eq:iser2}
\infty=1+\sum_{i=0}^{\infty}(p_i(y)-1) ,  \qquad \textrm{ if } f(y)\in\{\w,\co\}.
\end{equation}

\end{theorem}
\begin{proof}
We apply in turn Lemma \ref{series} to $f, f-p_0+1, f-(p_0+1)-(p_1+1),\ldots$
and obtain a sequence $p_0,p_1,p_2,\ldots $.

Put $$f_n=p_0+\sum_{i=1}^{n}(p_i-1).$$  

\smallskip
\noindent The inequality $p_{n+1}\leq p_n$ follows from the definition of $p_n$.

Fix $y\in[0,1]$ such that $f(y)$ is finite. For (1) it is enough to notice that
for some $m\in\w$
$$f_0(y)<f_1(y)<\ldots<f_{m-1}(y)<f_m(y)=f_{m+1}(y)=\ldots=f(y)$$

Fix $y\in[0,1]$ such that $f(y)=\infty$. For (2) it is enough to notice that
for each $n$, $p_n(y)\geq2$.

\medskip

\end{proof}

\medskip

%%%%%%%%%%%%%%%%%%%%%%%%%%%%%%%%%%%%%%%%%%%%%%%%%%%%%%%%%%%%%%%%%%%%%%%%%%%%%%%%
\bigskip

\noindent {\bf 3. The construction with the assumption $f^{\co}=\emptyset$}

We start with some preliminary notation.

\begin{figure}[h]%\label{picture}
\begin{center}
\scalebox{0.90}{\includegraphics{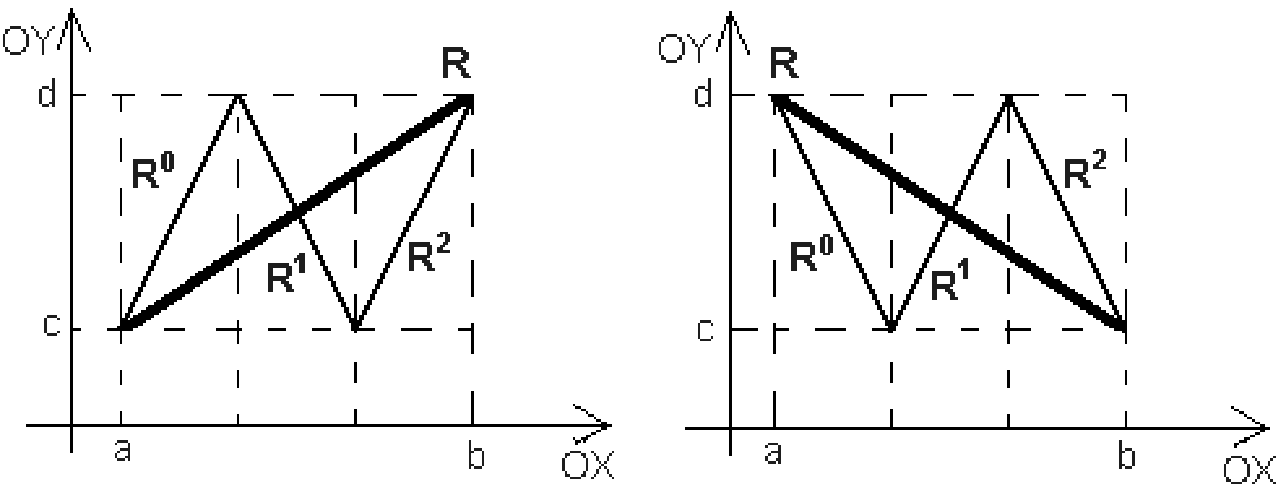}}

\end{center}
\caption{Graphs of $q_{\ra}$.}
	\label{fig:bloczek}

\end{figure}

Let $\ra$ be one of the diagonals in a rectangle $R=[a,b]\times[c,d]$.
For $e=0,1,2$, consider rectangles $[a+e\frac{b-a}{3},a+(e+1)\frac{b-a}{3}]\times[c,d]$ with 
  diagonals $\ra^e$ 
as shown in Figure \ref{fig:bloczek} and $q_{\ra}=\ra^0\cup \ra^1\cup \ra^2$.
Let $\pi_X$ and $\pi_Y$ be the projections onto the $x$ and $y$ axes.
Note that $\pi_X \ra=[a,b]$, $\pi_Y \ra=[c,d]$, $  \pi_Y \ra^0=[a,a+\frac{b-a}{3}]$, etc. 

For an open interval $I=(c',d')$ define
$$\ra\cdot I=\ra\cap ([a,b]\times[c',d']),$$
and for a family $P$ of open intervals, let
$$\ra\cdot P=\{\ra\cdot I\colon I\in P\}.$$
For a diagonal $\ra$ let
$$R=\pi_X \ra\times \pi_Y \ra$$
denote the rectangle from which $\ra$ arises. If $\Ra$ is a family of diagonals,
let $\RR$ be the corresponding family
$\{R\colon \ra\in\Ra\}$ of rectangles.

%RYSUNEK 1

We start with an idea of the construction. Fix a preindicatrix $f$.
To the sequence of simple preindicatrices 
$p_0,p_1,p_2,\ldots$ (obtained from a preindicatrix $f$ via Theorem \ref{series2})
 we will assign  families of open intervals 
$P_0,P_1,P_2,\ldots$, described below. Next, to each  interval  $I$ in $P=\bigcup_i P_i$ we assign
a graph  $q_{\ra}$  for some diagonal $\ra$ such that the closure of $I$ is equal to
$\pi_Y \ra$. All these parts glued together will
 form the required graph.

First, let us take any simple preindicatrix $p$. 
Let $\langle y_i\colon i\in\w \rangle$ be an enumeration of the endpoints of its 
intervals of type three. 
For every $y_i<y_j$, which are  endpoints of the same interval of type three:
\begin{itemize}
\item[$\bullet$]  if $p(y_i)=p(y_j)=1$, or if $p(y_i)=2$ and $p(y_j)=1$
      and $y_i$ is the endpoint of two intervals of type three, we choose an increasing sequence
      $(y_{i,k})_{k\in\mathbb{Z}}$ such that $\lim_{k\rightarrow\infty}y_{i,k}=y_j$ 
      and $\lim_{k\rightarrow -\infty}y_{i,k}=y_i$ 
      ($\mathbb{Z}$ denotes the set of integers).
 \item[$\bullet$]  if $p(y_i)=1$ and $p(y_j)=2$, or if $p(y_i)=p(y_j)=2$
      and $y_i$ is the endpoint of two intervals of type three, we choose an increasing sequence
      $(y_{i,k})_{k<0}$ such that $\lim_{k\rightarrow -\infty}y_{i,k}=y_i$ 
      and $y_{i,-1}=y_j$.
\item[$\bullet$] if $p(y_i)=2$ and $p(y_j)=1$
      and $y_i$ is the endpoint of exactly one interval of type three, we choose an increasing 
      sequence
      $(y_{i,k})_{k\geq 0}$ convergent to $y_j$, such that $y_{i,0}=y_i$.          
\item[$\bullet$] in the remaining cases put $y_{j,0}=y_i$ and $y_{j,1}=y_j$.
\end{itemize}

Let $P$ be the set of all intervals $(y_{n,k},y_{n,k+1})$. We assign 
in this manner $P_0$ to $p_0$, $P_1$ to $p_1$,...
Without loss of generality we assume that  each
interval in $P_{i+1}$ is has the length less then $\frac{1}{i+1}$, 
and is a subinterval of some interval in $P_i$ (if necessary, divide some interval
$(a,b)\in P_i$ into some intervals $(a,c_1), (c_1,c_2), (c_2,c_3),\ldots,(c_k,b)$
for some $k$).

\medskip

Now, assuming that $f^{\co}=\emptyset$, we describe the construction.
A graphical example illustrating the construction can be found after the description of the 
general construction at the end of Chapter 4.

\noindent {\emph{ Step 0.}}

We  put $\Ra_0=(id\restriction[0,1])\cdot P_0$

\noindent {\emph{Step n+1.}}
We have from the nth step  the family of diagonals $\Ra_n$.

 Define
$$\Ra_{n+1}=\bigcup_{\ra\in\Ra_n}\ra^0\cdot P_n.$$

We define now a sequence of functions $F_0,F_1,F_2\ldots$ mapping [0,1] onto [0,1].

 Let $F_0$ is $id\restriction[0,1]$, $F_1$ is $F_0$ modified in each $\ra\in\Ra_0$
to $q_{\ra}$, $F_{n+1}$ is $F_n$ modified in each $\ra\in\Ra_n$ to $q_{\ra}$ (note that we 
modify only the left part of each $\ra$). 
Finally we put $F=\lim_{n\to\infty}F_n$.

Note that $F_{n+1}=F_n$ outside $\bigcup_n\RR_n$ and that for each rectangle 
$R\in\RR_n$ the lenght 
of $\pi_Y R$ is less then $\frac{1}{n}$. 
This implies that the convergence is uniform and thereby $F$ is continuous.
$F$ maps [0,1] onto [0,1], $F(0)=0$ and $F(1)=1$. 

\medskip

Notice that each rectangle in $\RR_{n+1}$ is contained in some rectangle in $\RR_n$.
A decreasing sequence $R_0\supseteq R_1 \supseteq R_2\supseteq \ldots$ such that
for each $n$ $R_n\in\RR_n$ call a \textit{branch}. Note that $|\bigcap R_n|=1$.

Now fix $y\in[0,1]$. We shall prove that $F$ receives every value $y$ exactly $f(y)$
times. 

It is easy to prove by induction that $f_n$ defined as in the proof of  Theorem \ref{series2}
is the indicatrix of $F_{n+1}$.

If $f(y)=n$, then there exists such $m$, that for $k\geq m$ $f_k(y)=f(y)$.
Clearly if $F_m^{-1}(y)=\{a_1, a_2,\ldots,a_n\}$ then 
$F_k^{-1}(y)=F^{-1}(y)=\{a_1, a_2,\ldots,a_n\}$. Finally $\left| F^{-1}(y)\right| =n$.

Otherwise, if $f(y)=\w$, we can either find exactly one branch $(R_n)_n$ such that 
$\pi_Y[\cap R_n]=y$, or we can find exactly two branches $(R_n)_n$ and $(R_n')_n$  
such that $\pi_Y[\cap R_n]=\pi_Y[\cap R_n']=y$ (in this case for some $N\in\w$ we
 have $R_n=R_n'$
for $n<N$ and $R_n\neq R_n'$ for $n\geq N$; moreover, for $n\geq N$ $R_n$ are of the form 
$[a_n,b_n]\times [c_n,y]$ and $R_n'$ are of the form 
$[a_n',b_n']\times [y,d_n']$).

In the first case, for each $n$, there are  exactly two $x_n^a$, $x_n^b$ such that 
$(x_n^a,y), (x_n^b,y)\in R_n\setminus R_{n+1}$ and $F(x_n^a)=F(x_n^b)=y$.
In the second case, if $n<N$ as above we have exactly two $x_n^a$, $x_n^b$ such that 
$(x_n^a,y), (x_n^b,y)\in R_n\setminus R_{n+1}$ and $F(x_n^a)=F(x_n^b)=y$, and if
$n\geq N$ there is exactly one $x_n^a$ such that $(x_n^a,y)\in R_n\setminus R_{n+1}$ 
and $F(x_n^a)=y$ and there is exactly one $x_n^b$ such that
 $(x_n^b,y)\in R_n'\setminus R_{n+1}'$  and $F(x_n^b)=y$ (note that for $n>N$ we have
$x_n^a\neq x_n^b$).  

This observation and the above remark that branches have one-point intersection gives us 
$|F^{-1}(y)|=\w$.

\begin{remark}
If $rng f$ is bounded then $F=F_n$ for some $n\in\w$.
\end{remark}

%%%%%%%%%%%%%%%%%%%%%%%%%%%%%%%%%%%%%%%%%%%%%%%%%%%%%%%%%%%%%%%%%%%%%%%%%%%%%%%%%%%%%%%%%%%%%%
\medskip\bigskip

\noindent {\bf 4. The general construction }
%\enlargethispage{1000pt}

\bigskip

For a sequence $\alpha\in 3^\w$ write $W_\alpha$ for $\{n\in\w\colon \alpha(n)\neq 0\}$ and let
$$ \mathbb{P} =\{\alpha\in 3^{\w}\colon |W_\alpha|=\w \}.$$
The set $f^{\co}$ is  analytic, so if it is nonempty, there exists 
a continuous  mapping $\alpha\mapsto \alpha^*$ of $\mathbb{P}$ onto $f^{\co}$
such that $W_\alpha=W_\beta$ implies $\alpha^*=\beta^*$.
Such mapping can be obtained as a composition of a "projection" map of $\mathbb{P}$
onto $\mathbb{P}_2=\{\alpha\in 2^\w\colon \exists^{\infty}_n \alpha(n)\neq 0\} $
changing all 2's in sequences into 1's, and a mapping from $\mathbb{P}_2$ onto
$f^{\co}$ (note that $\mathbb{P}_2$ is a copy of the Baire space $\w^{\w}$). 
    
For a finite sequence $\tau\in3^{<\w}$ we put
$$\tau^*=\{\alpha^*\colon \tau\subseteq \alpha,\ \alpha\in \mathbb{P}\}.$$

We will also use the notation introduced in the previous chapter.

As for the case $f^{\co}=\emptyset$, to each  interval  in $P=\bigcup_i P_i$ we assign
a graph  $q_{\ra}$  for some diagonal $\ra$. All these parts glued together will
 form the required graph. 
We will use each of the intervals in $P$ exactly once.
If $I,J\in P$ and $I\subseteq J$, then J will be used before $I$. In general we will
not use at the step $i$ exactly intervals from $P_i$, $i\in\w$, as it was
 when $f^{\co}=\emptyset$.

Now we define families of diagonals $\ra_{\tau}$,
$\tau\in3^{<\omega}$ inductively along $|\tau|$. To do this, at the step $n+1$ of induction,
 we also define 
families of intervals $P_i^{n+1}, i\in\w$ and families of intervals
$P_i^{n,k}, k,i\in\w$ . Each diagonal in every $\ra_{\tau}$ will come from some (exactly one) 
 interval in $P$ and two different diagonals will come from different intervals in $P$.

\noindent {\emph{ Step 0.}}

We  put $\Ra_{\langle\rangle}=(id\restriction[0,1])\cdot P_0$ and we rename our families of
intervals as follows
$$P_0^0=P_1, P_1^0=P_2, \ldots, P^0_i=P_{i+1},\ldots$$

\noindent {\emph{Step n+1.}}
We have from the nth step:  the families of intervals $P_i^n$, $i\in\w$, 
and  the families of diagonals $\Ra_\tau$, $\tau\in 3^n$.

 Let $\langle \tau_k\colon k<3^n\rangle$ be an enumeration of $3^n$. 
Inductively along $k$ for $0\leq k<3^n$ we define families of diagonals $\ra_{\tau}$,
$|\tau|=n+1$, $\tau=\tau_k^\smallfrown e, e=0,1,2$; and
families of intervals
$P_i^{n,k}, i\in\w$ .

At the k-th step we are given (defined at the (k-1)-th step)  families of intervals $P_i^{n,k-1}$,
 $i\in\w$ (if $k=0$ we let $P_i^{n,-1}=P_i^n$, which are also given 
at the beginning of this step).

Define
$$\Ra_{\tau_{k}^\smallfrown 0}=\bigcup_{\ra\in\Ra_{\tau_{k}}}\ra^0\cdot P_0^{n,k}.$$
For $e=1,2$ define
$$\Ra_{\tau_{k}^\smallfrown e}=\bigcup\{\ra^e\cdot P_e^{n,k}\colon \ra\in\Ra_{\tau_{k}}, 
\tau_{k}^*\cap \pi_Y \ra\neq\emptyset\}. $$
Also,  to each $\ra\in\Ra_{\tau_{k}^\smallfrown e}$, $e=0,1,2$, assign a \textit{label} $\#R=e$.

Finally put for $i\in\w$ and all intervals $I$

 $I\in P_i^{n,k}$ \qquad iff

$ \bullet I\in P_{i}^{n,k-1} \textrm{ and there is no }  \ra\in\Ra_{\tau}
 \textrm{ such that } I\subseteq\pi_Y \ra$ , or 

$\bullet I\in P_{i+1}^{n,k-1} \textrm{ and  } I\subseteq\pi_Y \ra \textrm{ for some } \ra\in\Ra_{\tau}
\textrm{ such that } \tau^*\cap \pi_Y \ra=\emptyset$, or 

$\bullet I\in P_{i+3}^{n,k-1} \textrm{ and  } I\subseteq\pi_Y \ra \textrm{ for some } \ra\in\Ra_{\tau}
\textrm{ such that }  \tau^*\cap \pi_Y \ra\neq\emptyset.  $

\medskip

In other words: if $P_i'^{n,k-1}$ is $P_i^{n,k-1}$ without those intervals we have just used
in the step $(n,k)$ and if for some $i$ there is $I\in P_{i+1}'^{n,k-1}$ such that $I$ is not 
contained in (equivalently does not intersect)  any interval in $P_i'^{n,k-1}$,
 then we move $I$ 
to $P_i'^{n,k-1}$. We inductively iterate this operation. Finally, when for each 
$I\in P_{i+1}'^{n,k-1}$   
there is $J\in P_{i}'^{n,k-1}$ such that $I\subseteq J$, we put 
$P_i^{n,k}=P_i'^{n,k-1}$.

This finishes the step $k$.

At $k=3^n-1$,  let also $P_i^{n+1}=P_i^{n,k}$, $i\in\w $. 

This finishes  step $n+1$.

For $n\in\w$ let $$\RR_n=\bigcup_{\tau\in\ 3^n}\RR_\tau $$
 and $$\RR_{\w}=\bigcup_{n\in\w}\RR_n.$$

Note that $\RR_{\w}$ with inclusion becomes a tree. Accordingly, a decreasing
sequence of rectangles $(R_n)_n$, $R_n\in\RR_n$, 
will be called a \textit{branch} and the sequence of correspnding labels
$(\# R_{n+1})_n$ will be called a \textit{label} of this branch.
Note that for any branch $(R_n)_n$ 
there is $(x,y)$ such that $\{(x,y)\}=\bigcap_{n\in\w} R_n$.
We will say that $(R_n)_n$ \textit{converges} to $((x,y)$ $(R_n)_n\rightarrow (x,y)$).

\noindent{\bf Claim 1.} For any $\alpha\in 3^{\w}$ and $y\in[0,1]$ there are at most two
branches $ (R_n)_n$ that converge to $(x,y)$ for some $x$ and have label $\alpha$.
  
{\bf Proof: } Assume that there is a branch $(R_n)_n$ with a label $\alpha$ such that 
$(R_n)_n\rightarrow(x,y)$ for some $x$. Each $R_n$ is assigned to some interval $J_n$, 
namely $J_n=\pi_Y[R_n]$ (there is only one such assignment).
These intervals form a subsequence of some decreasing sequence $(I_k)_k$ such that 
$I_k\in P_k$ and $y\in I_k$, $k\in\w$. Note that there are at most two such sequences $(I_k)_k$ .
On the other hand, there is not more then one subsequence 
$(\widetilde{J}_n)_n$ of $(I_k)_k$ such that the sequence $(\widetilde{R}_n)_n$ 
of rectangles assigned to it form a branch of label $\alpha$. $\square$

\noindent{\bf Claim 2.} Suppose $(R_n)_n$ is a branch
with label $\alpha$.
Then if $W_\alpha$ is infinite then $(R_n)_n$ converges to ($x$,$\alpha^*$) for some $x$.

{\bf Proof: } For each $n\geq 1$ we have $R_n\in\RR_{\alpha|n}$ and
$R_{n+1}\in\RR_{\alpha|n^\sm \alpha(n)}$. 
If $n\in W_\alpha$ then $\alpha(n)\in\{1,2\}$ and therefore, because $R_{n+1}\subseteq R_n$,
$R_{n+1}$ witness that $(\alpha|n)^*\cap \pi_YR_n\neq\emptyset$. It
contains $\beta^*_n$ for some $\beta_n\supseteq \alpha|n$. 
Now $\lim_{n\in W_\alpha}\beta_n=\alpha$,
so by continuity $\lim_{n\in W_\alpha}\beta^*_n=\alpha^*$. Since also 
$\lim_{n\in W_\alpha}\beta_n^*\in\bigcap_{n\in W_\alpha}\pi_Y R_n=\bigcap_{n\in\w}\pi_Y R_n$,
then $\bigcap_{n\in\w}\pi_Y R_n=\{\alpha^*\}$ and
 we are done. $\square$

We define now a sequence of functions $F_0,F_1,F_2\ldots$ mapping [0,1] onto [0,1].
as follows:

 $F_0$ is $id\restriction[0,1]$, $F_1$ is $F_0$ modified in each $\ra\in\Ra_0$
to $q_R$, $F_{n+1}$ is $F_n$ modified in each $\ra\in\Ra_n$ to $q_R$
(an example of such $F_2$ can be found in Example 1). 
Finally we put $F=\lim_{n\to\infty}F_n$ and claim that this is our needed $F$.
The convergence is uniform, hence $F$ is continuous, maps [0,1] onto [0,1], 
and $F(0)=0$, $F(1)=1$. 

It remains to see that the sections of $F$ are as required.

Fix $y\in[0,1]$. Note that if $(R_n)_n\rightarrow (x,y)$ for some $x$ then $F(x)=y$.
Moreover, there are countably many $x$ such that $F(x)=y$ and $(x,y)$ is not a limit of any
branch. The reasoning is similar to that when we were proving that assuming $f^{\co}=\emptyset$
$f(a)=\w$ iff $|F^{-1}(a)|=\w$. Indeed, put $K_0=[0,1]^2\setminus\bigcup\RR_0$
and $K_{n+1}=\bigcup\RR_n\setminus\bigcup\RR_{n+1}$, $n\in\w$. Note that for every $n$
there are finitely many $x$ such that $(x,y)\in K_n$ and $F(x)=y$.

It is enough to prove: 

\begin{lemma}
Take $y\in[0,1]$.

\noindent (a) If $f(y)=n$ then there are exactly $n$ such $y$ that $F^{-1}(y)=n$.

\noindent (b) If $f(y)=\w$ then there are countably many $x$ such that $(x,y)$ is the 
limit of some branch.

\noindent (c) If $f(y)=\co$ then there is a perfect set $P$ such that 
for every $p\in P$ $(p,y)$ is the 
limit of some branch.
\end{lemma}
\begin{proof}
(a) As in the construction with $f^{\co}=\emptyset$.

(b) Let $(R_n)_n\rightarrow (x,y)$ and $\alpha$ be the label of this branch. By Claim 2
$W_\alpha$ is finite. There are countably many $\alpha$ such that $W_\alpha$ is finite.
Claim 1 finishes the proof.

(c) Fix $\alpha$ such that $y=\alpha^*$. Put $T=\{\beta\in3^{\w}\colon W_\beta=W_\alpha\}$
 and let 
$S=\{\beta\restriction n\colon \beta\in T, n\in\w\}$ be the tree corresponding to $T$.

We shall construct inductively a family $(R_\tau)_{\tau\in S}$ such that $R_\tau\in\RR_\tau$,
if $\tau\subseteq\sigma$ implies $R_\sigma\subseteq R_\tau$, $y\in\pi_Y R_\tau$ and
if $\beta\neq\beta'$ implies $\bigcap_n R_{\beta|n}\neq\bigcap_n R_{\beta'|n}$. This will 
finish the proof of (c).

First fix a sequence of intervals $(I_n)_n$ such that $I_n\in P_n$ and $y\in I_n$ 
for each $n$. Every $R_\tau$ will be assigned to some interval in this sequence.
As $R_\emptyset $ take the rectangle assigned to $I_0$ construted at the 0th step.
Assume that we already have $R_\tau$ for some $\tau$, $|\tau|=n$. If $\tau^{\smallfrown}0\in S$,
take $R_{ \tau^{\smallfrown}0 }\subseteq R_\tau$ constructed at the step (n+1).
If $\tau^{\smallfrown}1, \tau^{\smallfrown}2\in S$ take 
$R_{ \tau^{\smallfrown}1 }$, $R_{ \tau^{\smallfrown}2 }\subseteq R_\tau$ constructed at the 
step (n+1) (note that $\tau\in S$ implies $\tau^*\cap \pi_Y R_\tau\neq\emptyset$, 
which is witnessed by $y$). Note that it can happen that
$R_{\tau^{\smallfrown}1 }$ and $R_{\tau^{\smallfrown}2 }$ have common edge.
However, for all $e_1,e_2\in\{0,1,2\}$ we will have
$R_{ \tau^{\smallfrown}1^{\smallfrown}e_1 }\cap 
R_{ \tau^{\smallfrown}2^{\smallfrown}e_2 }\neq\emptyset$, whenever
$R_{ \tau^{\smallfrown}2^{\smallfrown}e_1 },
R_{ \tau^{\smallfrown}2^{\smallfrown}e_2 }\in S$.
In this manner we accomplish the last required condition on $(R_\tau)_\tau$.

\end{proof}

%%%%%%%%%%%%%%%%%%%%%%%%%%%%%%%%%%%%%%%%%%%%%%%%%%%%%%%%%%%%%%%%%%%%%%%%%%%%%%%%%%%%%%%%%%%%%%%%%%

%%%%%%%%%%%%%%%%%%%%%%%%%%%%%%%%%%%%%%%%%%%%%%%%%%%%%%%%%%%%%%%%%%%%%%%%%%%%%%%%%%%%%%%%%%%%
\bigskip
This ends our proof that conditions listed in Theorem \ref{main} are also sufficient.

\begin{example}
Let $f$ be a preindicatrix such that $f^{\co}\neq\emptyset$.

Let $p_0=3$ in (0,1) and $p_0(y)=2$ if $y\in\{0,1\}$.

Let $p_1=3$ in $(0,\frac{3}{4})$, $p_1(y)=2$ if $y\in\{0,\frac{3}{4}\}$
and $p_1=1$ in $(\frac{3}{4},1]$.

Let $p_2=3$ in $(0,\frac{1}{2})$, $p_2(0)=1$,  $p_2(\frac{1}{2})=2$ and
$p_2=1$ in $(\frac{1}{2},1]$.

Finally, let $p_3=3$ in $(\frac{1}{4},\frac{1}{2})$ and $p_3(y)=1$ 
for the remaining $y$.

The sketch of graph of $F_2$ is shown in  Figure \ref{fig:konstrukcja}

%RYSUNEK 2
\begin{figure}[h]%\label{picture}
\begin{center}
\scalebox{1}{\includegraphics{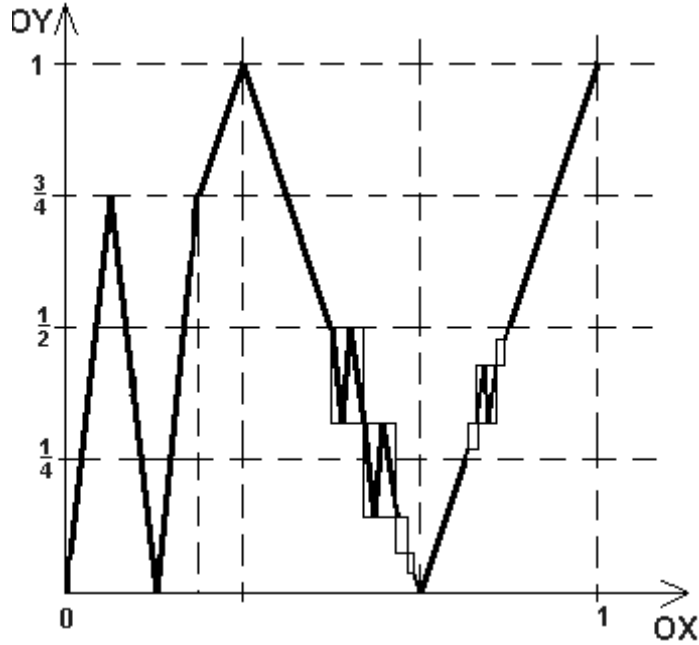}}

\end{center}
\caption{Sketch of graph of $F_2$.}
	\label{fig:konstrukcja}

\end{figure}

\end{example}

\bigskip
%%%%%%%%%%%%%%%%%%%%%%%%%%%%%%%%%%%%%%%%%%%%%%%%%%%%%%%%%%%%%%%%%%%%%%

\noindent {\bf 5. The general characterisation of indicatrices}

In the theorem below we present a general form of our result.

\begin{theorem}\label{main2}
  $f$ is an indicatrix  of a continuous function $F\colon[0,1]\rightarrow[0,1]$
 such that $F(0)=a$ and $F(1)=b$ for some $a\leq b$, if and only if:

\begin{enumerate}

\item If $f(y)\in\w$ and
 $y\in[0,a)\cup(b,1]$ then there exists $\varepsilon>0$ such that for any $y_1\in(y-\varepsilon,y)$
and $y_2\in(y,y+\varepsilon)$ we have $f(y_1)+f(y_2)\geq 2f(y)$
and if both of $f(y_1)$ and $f(y_2)$ are odd then 
$f(y_1)+f(y_2)>2f(y)$.

\item If $f(y)\in\w$ and
 $y\in(a,b)$ then there exists $\varepsilon>0$ such that for any $y_1\in(y-\varepsilon,y)$
and $y_2\in(y,y+\varepsilon)$ we have $f(y_1)+f(y_2)\geq 2f(y)$
and if both of $f(y_1)$ and $f(y_2)$ are even then 
$f(y_1)+f(y_2)>2f(y)$.

\item If $f(y)\in\w$ and
 $y=a=b,$ then there exists $\varepsilon>0$ such that for any $y_1\in(y-\varepsilon,y)$
and $y_2\in(y,y+\varepsilon)$ we have $f(y_1)+f(y_2)\geq 2(f(y)-1)$
and if both of $f(x)$ and $f(y)$ are odd then 
$f(y_1)+f(y_2)>2(f(y)-1).$

\item For all $y\in(0,1)$ such that $f(y)=\w$ and for all but countably many $y$  such that 
$f(y)=\co$,

\ for every  $n\in\w$ there exists $\varepsilon>0$ such that 
for any $y_1\in(y-\varepsilon,y)$
and $y_2\in(y,y+\varepsilon)$ we have 
$  f(y_1)+f(y_2)\geq n$. 

\item The set $f^{\co}$ is analytic. 

\end{enumerate}

If $y=0$ or $y=1$ we put in the above inequalities
 $f(y_1)=1$ and $f(y_2)=1$ respectively (and we also omit "$y_1\in(y-\varepsilon,y)$"
and "$y_2\in(y,y+\varepsilon)$" respectively).
\end{theorem}

%\enlargethispage{1000pt}

\begin{proof}
We start with the proof of necessity. Take $F$ as in the assumptions.
Consider $\widetilde{F}$ such that
$\widetilde{F}\restriction[0,\frac{1}{3}]$ and
$\widetilde{F}\restriction[\frac{2}{3},1]$ are linear functions linking
(0,0) with $(\frac{1}{3},a)$ and $(\frac{2}{3},b)$ with (1,1) respectively, and
its graph over $[\frac{1}{3},\frac{2}{3}]$ is a shrunk copy of $F$.
Let $\tilde{f}$ describes $\widetilde{F}$. $\tilde{f}$ satisfies the conditions
of Theorem \ref{main}. Now it can easily be checked that $f$ which describes $F$
satisfies the listed conditions.

Now we prove the sufficiency.
Consider $\hat{f}$ such that $\hat{f}(y)=f(y)-1$ for $y\in(0,a)\cup(b,1)$
and if $a\neq b$ also for $y\in\{a,b\}$, $\hat{f}(y)=f(y)-2$  if $y=a=b$, and
$\hat{f}(y)=f(y)$ for other $y$. Notice that $\hat{f}$ satisfies the conditions
of Theorem \ref{main}.  
Let $\widehat{F}$ such that
$\widehat{F}(0)=0$ and $\widehat{F}(1)=1$ be described by $\hat{f}$. 
Let $F$ be such that its graph over $[\frac{1}{3},\frac{2}{3}]$ 
is a shrunk copy of $\widehat{F}$, and $F\restriction[0,\frac{1}{3}]$
and $F\restriction[\frac{2}{3},1]$ are linear functions linking points 
$(0,a), (\frac{1}{3},0)$ and
$(\frac{2}{3},1),(1,b)$ respectively. $F$ is described by $f$.  

\end{proof}

Finally note that $f$ is an indicatrix  of  $F$
 such that $F(0)=a$ and $F(1)=b$ for some $a\geq b$ if and only if it is 
an indicatrix  of  $G$ such that   for $x\in[0,1]$ we have $G(x)=F(1-x)$.
Such a $G$ satisfies $G(0)\leq G(1)$.

\bigskip

{\bf Acknowledgments.}
I would like to thank E. Damek and J. Pawlikowski for their help in improving the exposition
of the paper and to Szymon Glab for many remarks concerning the exposition. 
 I also would like to thank M. Balcerzak, 
for a valuable review of an earlier version of the paper.

\end{document}